\documentclass[a4paper,10pt]{article}

\usepackage{amsfonts,amssymb,mathrsfs,amscd}
\def \R {{\mathbb {R}}}

\def \Z {{\mathbb {Z}}}

\def\eps{\varepsilon}

\usepackage[T1]{fontenc} % иногда T1 вместо T2A
\usepackage[cp1251]{inputenc}  
\usepackage[russian]{babel}

\title{\bf Around  Krygin-Atkinson, Shneiberg theorems, \\ the recurrence with  zero integrals}
\author{Valery V. Ryzhikov}
\date{}
\usepackage{graphicx}
\graphicspath{{pictures/}}
\DeclareGraphicsExtensions{.pdf,.png,.jpg}
\begin{document}
\large
\maketitle
%\begin{abstract}  
%Theorems related to the old result of Krygin and Atkinson
%on conservative cylindrical cascades are discussed.
%Библиография: 15 названий, УДК 517.987  

%Ключевые слова и фразы: 
%\end{abstract} 

\bf Abstract. \rm We recall  theorems by  Krygin, Atkinson, Shneiberg and propose the following assertion.  Let $T_t$ be an  ergodic flow on $(X,\mu)$, let a function $f:X\to\R$ have zero mean, and $\mu(A)>0$ for a set $A\subset X$. Then for almost all $x\in A$ with $f(x)\neq 0$ there exists a sequence $t_k\to\infty$ such that $\int_0^{t_k}f(T_sx)ds=0$ and $T_{t_k}x\in A$.
\rm
\section{Introduction}
Let $S$ be a probability measure-preserving invertible transformation (automorphism)
of the space $(X, \mu)$, $f:X\to\Z$ be a $\mu$-measurable function.
The cylindrical cascade $C: X\times \Z\to X\times \Z$ is defined by the formula
$$ C(x,z)=\left(Sx, z+f(x)\right),$$
it preserves the measure $\bar \mu =\mu\times\sharp$.

Let $T_t$ be a flow preserving the measure $\mu$. The cylindrical flow $C_t: X\times \R\to X\times \R$ is defined similarly:
$$ F_t(x,r)=\left(T_tx, r+ \sigma(t,x)\right),$$
$$\sigma(t,x):=\int_{0}^{t}f(T_sx)\,ds.$$
It preserves the measure $\bar \mu =\mu\times m$, where $m$ is the Lebesgue measure on $\R$.

Krygin and Atkinson \cite{Kr},\cite{At} independently proved that the cylindrical cascade with ergodic $S$ and integrable function $f:X\to\Z$ with zero mean
are recurrent. This means that for any set
$B\subset A$ of positive measure, for almost all points $x\in B$,
the points $(x,z)$ return to $B\times \{z\}$ infinitely many times under the action of the cylindrical cascade. Since
$$ C^n(x,z)=\left(Sx, z+\sum_{i=0}^{n-1}f(T^ix)\right),$$
for Birkhoff sums, $$\sum_{i=0}^{n_k-1}f(T^ix)=0$$
for an infinite set of times $n_k$.
Krygin formulated his result in the case
when $S$ is an ergodic rotation, but in the proof he used only the ergodicity property of this automorphism. Atkinson noted this in the abstract of Zbl 0342.60049. In the case where $A$ is additionally a metric space and every open set has positive measure, then the Krygin-Atkinson theorem directly implies for almost all $x\in X$ the existence of a sequence $n_k\to\infty$ for which $T^{n_k}x\to x$ and
$\sum_{i=0}^{n-1}f(T^ix)=0$ hold.

%\newpage
The case of flows is considered in the paper by Shneiberg \cite{Sh}, in particular,
he proved that for an ergodic flow $T_t: X\to X$ and $f:X\to\R$ with zero mean for almost all $x\in X$ there exists a
sequence $t_k\to\infty$
$$ \sigma(t_k,x)=\int_0^{t_k} f(T_sx)\, ds =0.$$

In the paper by Denisova \cite{De}
the question of when $\sigma(t,x)$ has infinitely many zeros
$t_k\to\infty$ in each neighborhood of the point $x$ is studied.

\vspace{3mm}
\bf Theorem (\cite{De}). \it Let $T_t$ be an ergodic flow on a compact metric space $X$ with finite Caratheodory measure $\mu$, and let
$f:X\to\R$ be continuous and have zero mean. Then for almost all $x$ such that $f(x)\neq 0$ there exists a sequence $t_k\to\infty$ such that $\sigma(t_k,x)=0$ and $T_{t_k}x\to x$.\rm

\vspace{3mm}
\bf Theorem A. \it Let $T_t$ be a special ergodic flow on $(X,\mu)$, let $f:X\to\R$ have zero mean, $\mu(A)>0$. Then for almost all $x\in A$ with $f(x)\neq 0$ there exists a sequence $t_k\to\infty$ such that $\sigma(t_k,x)=0$ and 
$T_{t_k}x\in A$.\rm

\vspace{3mm}
\bf  Corollary. \it Let $T_t$ be an ergodic flow with respect to the measure $\mu$  on a metric space $X$, and let the measure $\mu$ of each open set be positive. Suppose that for a function $f:X\to\R$ with zero mean for almost all $x$ the integrals $\sigma(t,x)$ make sense and are continuous in $t$. Then for almost all $x$ with $f(x)\neq 0$ there is a sequence $t_k\to\infty$ such that $\sigma(t_k,x)=0$ and $T_{t_k}x\to x$.\rm

\section{Recurrence with zero integrals}
A measurable measure-preserving flow is isomorphic to a special flow $T_t$.
We consider the case when the orbits have zero measure.
The phase space $X$ for $T_t$ is the part of the plane under the graph of an integrable non-negative function $r: A\to \R^+$, $A\subset \R$. The metric on the phase space is the metric $\rho$ on the plane, and it is locally preserved by the flow.
Recall that the point $x=(a,b)$ of the phase space moves vertically upward with a constant velocity: $T_t(a,b)=(a, b+t)$. The boundary of the phase space is glued in such a way that the points $(a, r(a))$ are identified with the points $(P(a),0)$, where $P$ is a given automorphism on $A$, the measurable domain of $r$.

\vspace{3mm}
\bf Theorem B. \it Let $T_t$ be a special ergodic flow on $(X,\mu)$ (or ergodic winding of a torus).  For a function $f:X\to\R$ with zero mean for almost all $x$, for $f(x)\neq 0$ there exists a sequence $t_k\to\infty$ such that $\sigma(t_k,x)=0$ and $T_{t_k}x\to x$.\rm

\vspace{3mm}
Proof. The function $f$ is integrable on $X$, from Fubini's theorem it follows that the functions $F(s,x)= f(S_sx)$ for almost all $x$, as functions of $s$, are integrable on vertical segments with endpoints
$(a,0)$ and $(a,r(a))$. Therefore, for almost all $x$, the indefinite integral
$$\Phi(t)= \int_0^{t} f(T_sx)\, ds$$ is an absolutely continuous function
on $\R$, hence for almost all $t$ there exists a limit
$$ \lim_{\Delta \to 0}\frac {\int_t^{t+\Delta}f(T_sx)\, ds}{\Delta}=f(T_tx).$$
But then for almost all
$x$ there exists a limit
$$ \lim_{ t\to 0}\frac {\int_0^{t}f(T_sx)\, ds}{t}=f(x).$$
It follows from this that for almost all $x$ with $f(x)\neq 0$ there exists a number $\delta(x)>0$ such that $$ \int_0^{\delta(x)}f(T_sx)\, ds >f(T_tx)\delta(x)/2.$$
We can obviously choose the function $\delta(x)$ in a measurable way, where $\delta(x)<\delta$
for a predetermined $\delta>0$. Let $U$ denote the intersection of the small disk with the set of $x$-s for which we have defined $\delta(x)$. Let the measure $U$ be positive.
In the phase space of a cylindrical flow, we consider the set
$$ E=\{(x,h)\,:\, x\in U, f(x)>0, 0>h>-\delta(x)/4\}.$$
Since the cilindrical  flow $F_t$  is recurrent, for any $N$ there is $t>N$
$$ \bar\mu (F_tE\cap E)>0.$$ Let
$$(T_tx,-h)\in F_tE\cap E, \ \  0<h<\delta(x)/4,$$     
then 
%for all $x\in T_{-t}U\cap U$ it is true that $T_t x\in U$ and
 by the continuity of the integral there is $\Delta(x,h)\in (0,\delta(x))$ such that
$$ \int_0^{\Delta(x,h)}f(T_s T_tx)\, ds=h.$$
This means that
$$\int_0^{t+\Delta(x,h)}f(T_s T_tx)\, ds=0,$$
and by choice of $U$  the points $x$ and $T_{t+\Delta(x,h)}x$ may be as close as we want.
Let's summarize the above in the following form.

\vspace{3mm}
\bf Lemma. \it In each set $A\subset \{x: f(x)> 0\}$ of positive measure for every $N$   we find positive measure set of  $x\in A$ such that  there is  $t_N(x)>N$ for which  the distance between
$x$ and $T_{t_N(x)}x$ is less than $1/N$ and
$$ \sigma(t_{N(x)},x):=\int_0^{t_{N(x)}}f(T_s x)\, ds=0.$$ \rm

\vspace{3mm}
We call good a points $x$ for which $f(x)\neq 0$ and there is a sequence $t_N\to\infty$ for which $\sigma(t_N,x)=0$ and $\rho(T^{t_N}x,x)<1/N$.
A point $x$ is bad if there exists $N$ such that it is $N$-bad: for all $t>N$
from $\rho(T^{t_N}x,x)<1/N$ it follows that $\sigma(t_N,x)\neq 0$.
If the measure of bad points is positive, then the measure of $N$-bad points is positive for some  $N$. But by the lemma there is a point in it that is not $N$-bad. Thus, the set of bad points has measure 0.

The theorem B is proved for points where $f(x)>0$, to finis the proof we reason symmetrically for $f(x)<0$.
For an  ergodic winding of the torus as $T_t$, we get the same with the same  proof.

\vspace{3mm} 
\bf Proof of Theorem A. \rm The complete proof was given in \cite{11}, here we will briefly outline the main points. We call a point $x\in A$, $\mu(A)>0$, good for $A$ with respect to the ergodic flow $T_t$ and a nonzero on $A$  function $f$ with zero mean, if in its orbit for every $N$ there exist points $T_{t_1}x\in A$, $T_{t_1+t}x\in A$ such that
$$\int_{t_1}^{t_1+t}f(T_s x)ds =0, \ \ t_1\geq 0, \ t>N.$$
Using, for example, Theorem B and the properties of absolutely continuous functions,  one can show that almost all points from $A$ are good. But then
$\mu(A_N)=0$ for  the set $A_N$ of those points  $x\in A$ for which, for all  $t>N$ such that $T_tx\in A$ we have  $\int_{0}^{t}f(T_s x)ds\neq 0$. Indeed,
if $\mu(A_N)>0$, then almost all of its points are good for $A_N$, and this, as is easy to see, will contradict the definition of $A_N$.

\section{Remarks}
\bf  Infinite  Krygin-Atkinson Theorem. \rm The following statement is easy to prove by reducing the infinite case to the case of finite measure.

\vspace{3mm}
\bf Theorem C. \it Let $(X,\mu)$ be a standard space with sigma-finite measure, and $S$ an ergodic automorphism of this space. If $f:X\to\Z$ has zero mean, then
the cylindrical cascade $C: X\times \Z\to X\times \Z$,
$ C(x,n)=\left(Sx, n+f(x)\right),$
is conservative. \rm

\vspace{3mm}
Consider the induced cylindrical cascade on the set $A\times \Z$,
$\mu(A)=1$, $A\subset X$.
Let $x\in A$,  $n(x)$ satisfy $S^{n(x)}\in A$, $S^k\notin A,$ for $ 0<k<n(x)$.
Put
$$\tilde Sx= S^{n(x)}x, \ \ \ \tilde f (x)=\sum_{i=0}^{n(x)-1}f(S^ix),
\ x\in A.$$
We have: $\tilde S:A\to A$ is ergodic on $A$ and $\int_A \tilde f (x)\, d\mu =0$.
The induced cylindrical cascade $\tilde C:A\times\Z\to A\times \Z$ is recurrent,
therefore so is the original cascade $C$.

\vspace{3mm}
The bibliography \cite{De} lists many works related to the topic of recurrence of cylindrical cascades, see \cite{Mo}--\cite{Ko}.
In addition, we note that the Krygin-Atkinson theorem was used by the author to lift the multiple mixing property and to lift the triviality property of pairwise independent joinings, called the del Junco-Rudolph property, see \cite{23}. This property plays an important role in the study of Rokhlin's multiple mixing problem.

Let $n>2$, an $S^{\otimes n}$-invariant measure on the cube $X^{\times n}$ with projections $\mu^{\otimes 2}$  onto two-dimensional faces.
The JR property for $S$ is that each such measure on $X^{\times n}$ is trivial, i.e. equal to $\mu^{\otimes n}$.

\vspace{2mm}
{\bf Theorem (\cite{97},\cite{23}). \it Let $R:X\times Y\to X\times Y$ be the skew product over $S$ defined by
$$
R(x,y)=(S(x),T^{n(x)}(y)), \quad \int n(x)d\mu = 0.
$$
If $R,T$ have mixing and $S$ mixes with multiplicity $k$, then the skew product $R$ inherits mixing of multiplicity $k$.

If $R,T$ mix and $S$ has the JR-property, then $R$ inherits the JR-property.\rm

%\newpage
We also note, following Benjamin Weiss, that the integrability of $f$ is not necessary for the recurrence of the cylindrical cascade. Theorem 1.4 \cite{We} implies the following assertion.

\vspace{2mm}
\bf Theorem D. \it Let $(X,\mu)$ be a probability space and $S$ an automorphism of this space. If $f:X\to\Z$ is such that
$$\mu\left( x: \left|\sum_0^{n-1} f(T^ix)\right|>\eps n\right)\to 0, \ n\to\infty,$$
with probability 1 an infinite number of times we observe $\sum_0^{n-1} f(T^ix)=0$. \rm

\large

\newpage 
{\bf \LARGE   Вокруг теорем  Крыгина-Актинсона и Шнейберга,\\
возвратность траекторий потока с нулевым интегралом.}  

\vspace{3mm} 
Пусть $S$ -- сохраняющее вероятностную меру обратимое преобразование (автоморфизм)
 пространства $(X, \mu)$, $f:X\to\Z$ -- $\mu$-измеримая функция.  

Цилиндрический каскад  $C: X\times \Z\to X\times \Z$ определен формулой
$$ C(x,z)=\left(Sx, z+f(x)\right),$$
сохраняет меру $\bar \mu =\mu\times\sharp$.
 
Пусть $T_t$ -- поток, сохраняющий меру $\mu$.  Цилиндрический поток  $C_t: X\times \R\to X\times \R$ задается аналогично:
$$ F_t(x,r)=\left(T_tx, r+ \sigma(t,x)\right),$$
$$\sigma(t,x):=\int_{0}^{t}f(T_sx)\,ds.$$ 
Он сохраняет меру  $\bar \mu =\mu\times m$,  где $m$ -- мера Лебега на $\R$.

Крыгин и Аткинсон \cite{Kr},\cite{At} независимо доказали, что цилиндрический каскад с эргодическим $S$ и интегрируемой функцией $f:X\to\Z$ с нулевым средним
являются рекуррентными.  Это означает, что для всякого множества  
$B\subset A$ положительной меры для почти всех точек $x\in B$ верно, что 
точки $(x,z)$ под действием цилиндрического каскада бесконечное число раз
возвращаются в $B\times \{z\}$.  Так как 
$$ C^n(x,z)=\left(Sx, z+\sum_{i=0}^{n-1}f(T^ix)\right),$$
 для сумм Биркгофа выполнено  $$\sum_{i=0}^{n_k-1}f(T^ix)=0$$
 для бесконечного множества моментов времени $n_k$.  
   Крыгин сформулировал свой результат в случае,
когда $S$ -- эргодический поворот, но в доказательстве использовал только свойство эргодичности этого автоморфизма. Аткинсон отметил это в автореферате Zbl 0342.60049. В случае, когда $A$ дополнительно является метрическим пространством и всякое открытое множество имеет положительную меру, то непосредственно из теоремы 
Крыгина-Аткинсона вытекает для почти всех  $x\in X$ существование   последовательности $n_k\to\infty$, для которой $T^{n_k}x\to x$ и 
$\sum_{i=0}^{n-1}f(T^ix)=0$.

%\newpage
Случай потоков рассмотрен   в работе Шнейберга \cite{Sh}, в частности,
им доказано, что для эргодического потока $T_t: X\to X$  и $f:X\to\R$ с нулевым средним  для почти всех $x\in X$  найдется 
последовательность $t_k\to\infty$
$$ \sigma(t_k,x)=\int_0^{t_k} f(T_sx)\, ds =0.$$

%\newpage
В статье Денисовой \cite{De}  
изучается вопрос о том, когда $\sigma(t,x)$ имеет бесконечно много нулей 
$t_k\to\infty$  в каждой окрестности  точки $x$.

\vspace{3mm}
\bf Теорема  (\cite{De}). \it Пусть  $T_t$ -- эргодический поток на компактном  метрическом пространстве $X$  с конечной мерой Каратеодори $\mu$, а  функция
$f:X\to\R$  непрерывна и имеет  нулевое среднее значение. Тогда для почти всех $x$ таких, что  $f(x)\neq 0$  найдется  последовательность $t_k\to\infty$ такая, что $\sigma(t_k,x)=0$ и $T_{t_k}x\to x$.\rm

\vspace{3mm}
\bf Теорема A. \it Пусть  $T_t$ -- специальный эргодический поток на $(X,\mu)$,    функция $f:X\to\R$  имеет    нулевое  среднее, $\mu(A)>0$. Тогда для почти всех $x\in A$ при  $f(x)\neq 0$   найдется  последовательность $t_k\to\infty$, для которой $\sigma(t_k,x)=0$ и $T_{t_k}x\in A$.\rm

\vspace{3mm}
\bf Следствие. \it  Пусть  $T_t$ -- эргодический относительно меры $\mu$ 
поток на   метрическом пространстве $X$, мера $\mu$ каждого открытого множества положительна. Предположим, что для   функции $f:X\to\R$  с   нулевым  средним   для почти всех $x$ интегралы $\sigma(t,x)$
имеют смысл и непрерывны по $t$. Тогда для почти всех $x$ при  $f(x)\neq 0$   найдется  последовательность $t_k\to\infty$, для которой $\sigma(t_k,x)=0$ и $T_{t_k}x\to x$.\rm

\vspace{3mm}
\bf Возвратность точек и нулей интегралов вдоль траекторий. \rm
 Измеримый поток, сохраняющий меру,  изоморфен специальному потоку $T_t$.
Мы рассматриваем случай, когда орбиты имеют меру 0.  
Фазовое пространство  $X$ для  $T_t$  является частью плоскости под графиком интегрируемой неотрицательной функции $r: A\to \R^+$, $A\subset \R$. Метрика на фазовом пространстве есть метрика $\rho$  на плоскости и она локально сохраняется потоком.
Напомним, что точка $x=(a,b)$ фазового пространства движется вертикально вверх с постоянной скоростью: $T_t(a,b)=(a, b+t)$. Граница фазового пространства склеена таким образом, что точки $(a, r(a))$ отождествляется с точками $(P(a),0)$, где $P$ -- заданный автоморфизм на $A$ -- измеримой области определения функции $r$.

\vspace{3mm}
\bf Теорема B. \it Пусть  $T_t$ -- специальный эргодический поток на $(X,\mu)$ (или эргодическая обмотка тора). 
 Для   функции $f:X\to\R$  с   нулевым  средним   для почти всех $x$ при  $f(x)\neq 0$   найдется  последовательность $t_k\to\infty$, для которой $\sigma(t_k,x)=0$ и $T_{t_k}x\to x$.  \rm

\vspace{3mm}
Доказательство.  Функция $f$ интегрируема на $X$,    bз теоремы Фубини следует, что функции  $F(s,x)= f(S_sx)$ для почти всех $x$
как функции от $s$   интегрируемы на вертикальных отрезках с концами
$(a,0)$ и $(a,r(a))$.
Поэтому для почти всех $x$ неопределенный интеграл  
$$\Phi(t)= \int_0^{t} f(T_sx)\, ds$$ является  абсолютно непрерывной функцией
на $\R$, следовательно для почти всех  $t$  существует предел
$$ \lim_{\Delta t\to 0}\frac  {\int_t^{t+\Delta}f(T_sx)\, ds}{\Delta}=f(T_tx).$$
Но тогда для почти  всех 
  $x$  существует предел
$$ \lim_{ t\to 0}\frac  {\int_0^{t}f(T_sx)\, ds}{t}=f(x).$$
Из этого вытекает, что для почти всех $x$ при $f(x)\neq 0$ найдется число $\delta(x)>0$ такое,
что $$ \int_0^{\delta(x)}f(T_sx)\, ds >f(T_tx)\delta(x)/2.$$
Функцию $\delta(x)$ мы, очевидно, можем выбрать измеримым образом, причем $\delta(x)<\delta$
для наперед заданного $\delta>0$.  Пусть $U$ обозначает  пересечение маленького диска с множеством таких $x$-ов, для которых мы определили $\delta(x)$. Пусть мера $U$ положительна.
В фазовом пространстве цилиндрического потока расмотрим множество
$$    E=\{(x,h)\,:\,   x\in U, f(x)>0,   0>h>-\delta(x)/4\}.$$
Так как цилиндрический  поток $F_t$  рекуррентный, для всякого $N$ найдется  найдется  $t>N$ такое, что 
$$ \bar\mu (F_tE\cap E)>0.$$ Пусть
$$(T_tx,-h)\in F_tE\cap E, \ \  0<h<\delta(x)/4,$$  
%Тогда для всех $x\in  T_{-t}U\cap U$ верно, что $T_t x\in U$ и 
в силу непрерывности интеграла  найдется $\Delta(x,h)<\delta(x)$ такое, что 
$$  \int_0^{\Delta(x,h)}f(T_s T_tx)\, ds=h.$$
Это означает, что 
 $$\int_0^{t+\Delta(x,h)}f(T_s T_tx)\, ds=0,$$
причем $x$ и $T_{t+\Delta(x,h)}$ близки.
Подведем итог сказанному.   

\vspace{3mm}
\bf Лемма. \it  В каждом множестве $A\subset \{x: f(x)> 0\}$ положительной меры  для всякого $N$ и  найдется множество положительной меры таких  $x$, что для некоторого  $t_N(x)>N$ расстояние между 
$x$ и $T_{t_N(x)}x$ меньше $1/N$ и 
$$ \sigma(t_{N(x)},x):=\int_0^{t_{N(x)}}f(T_s x)\, ds=0.$$  \rm 

\vspace{3mm}
Назовем хорошими  точки $x$, для которых    $f(x)\neq 0$ и   найдется  последовательность $t_N\to\infty$, для которой $\sigma(t_N,x)=0$ и $\rho(T^{t_N}x,x)<1/N$.
Точка $x$ плохая, если найдется такое $N$, что она  $N$-плохая:  для всех $t>N$ 
из $\rho(T^{t_N}x,x)<1/N$  вытекает $\sigma(t_N,x)\neq 0$.
Если мера плохих точек положительна, то мера $N$-плохих точек положительна для некоторого 
$N$. Но в силу леммы в ней найдется точка, которая не является $N$-плохой. Таким образом, множество плохих точек имеет меру 0.

Теорема доказана для точек, где $f(x)>0$, симметрично рассуждаем при $f(x)<0$. Теорема B доказана.

Если в качестве $T_t$ взять эргодическую обмотку тора, то получаем
идентичный результат с идентичным доказательством.

\bf Доказательство теоремы А. \rm Полное    доказательство  было дано в \cite{11}, здесь мы изложим кратко основную идею.   Назовем точку $x\in A$,  $\mu(A)>0$, хорошей для $A$ относительно эргодического потока $T_t$ и функции $f$ с нулевым средним  и  ненулевой на $A$, если   в ее орбите для всякого $N$ найдутся точки  $T_{t_1}x\in A$,  $T_{t_1+t}x\in A$ такие, что 
$$\int_{t_1}^{t_1+t}f(T_s x)ds =0,  \ \ t_1\geq 0, \ t>N.$$
Пользуясь, например,  теоремой В и свойствами абсолютно непрерывных функций,
можно показать, что почти все точки  множества $A$ хорошие. Но тогда
множество $A_N$, состоящее из  тех точек $x$ из $А$, для которых при $T_tx\in A$  и $t\geq N$ всегда выполнено 
$$\int_{0}^{t}f(T_s x)ds\neq 0,$$   имеет меру 0. В самом деле, 
если  $\mu(A_N)>0$, то почти все его точки хорошие для $A_N$, а это, как несложно заметить,  будет противоречить  определению $A_N$.

\vspace{5mm}
\bf Бесконечная теорема Крыгина-Аткинсона.  \rm 

\vspace{3mm}
\bf Теорема C.  \it  Пусть $(X,\mu)$ -- стандартное пространство с сигма-конечной мерой, $S$ -- эргодический автоморфизм этого пространства.  Если $f:X\to\Z$  имеет нулевое среднее, то 
цилиндрический каскад  $C: X\times \Z\to X\times \Z$,
$ C(x,n)=\left(Sx, n+f(x)\right),$ 
является консервативным. \rm

\vspace{3mm} 
Это  утверждение несложно доказать, сведя  случай бесконечной меры 
к случаю конечной меры.  Рассмотрим индуцированный цилиндрический каскад на множестве $A\times \Z$,
$\mu(A)=1$, $A\subset X$. 
Пусть $x\in A$ и для $n(x)$ выполнено $S^{n(x)}\in A$, $S^k\notin A,$ при $ 0<k<n(x)$.
Положим 
$$\tilde Sx= S^{n(x)}x, \ \ \ \tilde f (x)=\sum_{i=0}^{n(x)-1}f(S^ix),
\ x\in A.$$ 
Имеем: $\tilde S:A\to A$  эргодичен на $A$ и $\int_A  \tilde f (x)\, d\mu =0$.
Индуцированный цилиндрический каскад $\tilde C:A\times\Z\to A\times \Z$ рекуррентен,
следовательно, таков и исходный каскад $C$.

\vspace{3mm} 
В библиографии  \cite{De} указан ряд работ,  связанных с тематикой рекуррентности цилиндрических каскадов, см. \cite{Mo}--\cite{Ko}. 
В дополнение к этому, отметим, что теорему Крыгина-Аткинсона автор использовал для поднятия свойства кратного перемешивания  и поднятия   свойства тривиальности джойнингов с попарной независимостью, названного свойством  дель Джунко-Рудольфа, см. \cite{23}.  Это свойство играет важную роль в исследовании проблемы Рохлина о кратном перемешивании.
Пусть $n>2$, $S^{\otimes n}$-инвариантная мера на кубе с проекциями $\mu^{\otimes 2}$ на двумерные грани называется попарно независимым самосоединением автоморфизма $S$.
JR-свойство  для автоморфизма $S$ означает, что каждое такое самоприсоединение  на $X^{\times n}$ тривиально, является мерой  $\mu^{\otimes n}$.

\vspace{2mm}
{\bf Теорема (\cite{97}).  \it Пусть $R:X\times Y\to X\times Y$ -- косое произведение над автоморфизмом $S$, заданное формулой
$$
  R(x,y)=(S(x),T^{n(x)}(y)), \quad \int n(x)d\mu = 0.
$$
Если $R,T$ обладают перемешиванием, а автоморфизм  $S$ перемешивает с кратностью $k$, то  косое произведение $R$ наследует  перемешивание  кратности $k$.

Если $R,T$ перемешмвают и $S$ обладает   JR-свойством, то $R$ наследует JR-свойство.\rm

\vspace{2mm}
Отметим также, следуя Бенжамину Вейсу, что  интегрируемость функции $f$ не обязательна для рекуррентности цилиндрического каскада.   Из теорема 1.4  \cite{We} вытекает следующее утверждение.  

\vspace{2mm}
 \bf Теорема D.  \it  Пусть $(X,\mu)$ -- вероятностное  пространство, $S$ -- автоморфизм этого пространства.  Если $f:X\to\Z$  такова, что 
$$\mu\left( x: \left|\sum_0^{n-1} f(T^ix)\right|>\eps n\right)\to 0, \ n\to\infty,$$ 
 с вероятностью 1  бесконечное число раз мы наблюдаем $\sum_0^{n-1} f(T^ix)=0$.  \rm


\normalsize
\begin{thebibliography}{99}

\bibitem{Kr} 
A.B. Krygin, An example of a cylindrical cascade with anomalous metric properties, Moscow State University Bulletin, ser. 1, 1975, no. 5, 26-32;


\bibitem{At} G. Atkinson, Reccurence of co-cycles and random walks, J. London Math. Soc., 13 (1976), 486-488 


\bibitem{Sh} 
I.Ya. Shneiberg, Zeros of integrals along trajectories of ergodic systems, Funct. Anal. Appl., 19:2 (1985), 160-161;

\bibitem{De}
N. V. Denisova, Recurrence of integrals of conditionally periodic functions, Dokl. Math., 108:1 (2023), 316-319;

\bibitem{11} V.V. Ryzhikov, Recurrence of integral zeros on trajectories of ergodic flow. Preprint (2024)


\bibitem{Mo}
N.G. Moshchevitin, Recurrence of the integral of a smooth three-frequency conditionally periodic function, Math. Notes, 58:5 (1995), 1187-1196

\bibitem{Kon}
S.V. Konyagin, Recurrence of the integral of an odd conditionally periodic function, Math. Notes, 61:4 (1997), 473-479


\bibitem{Koz}
V. V. Kozlov, Methods of qualitative analysis in the dynamics of a rigid body.
Scientific Publishing Center Regular and Chaotic Dynamics, Izhevsk, 2000. 248 pp. 

\bibitem{Ko}
A.V. Kochergin, On the Growth of Birkhoff Sums over a Rotation of the Circle, Math. Notes, 113:6 (2023), 784-793


\bibitem{97}
V.V. Ryzhikov, 
Polymorphisms, joinings, and the tensor simplicity of dynamical 
systems, Funct. Anal. Appl., 31:2 (1997), 109-118

\bibitem{23}
V.V. Ryzhikov, Self-joinings and generic extensions of ergodic systems, Funct. Anal. Appl., 57:3 (2023), 236-247 

\bibitem{We} B. Weiss, 
Single orbit dynamics. BMS Regional Conference Series in Mathematics 95. Providence, RI: AMS, 113 p. (2000).

\end{thebibliography}

\normalsize
\begin{thebibliography}{99}

\bibitem{Kr} 
A.B. Krygin, An example of a cylindrical cascade with anomalous metric properties, Moscow State University Bulletin, ser. 1, 1975, no. 5, 26-32;

А.Б. Крыгин, Пример цилиндрического каскада с аномальными метрическими свойствами, Вестник МГУ, сер. 1, 1975, № 5, 26-32

\bibitem{At} G. Atkinson, Reccurence of co-cycles and random walks, J. London Math. Soc., 13 (1976), 486-488 


\bibitem{Sh} 
I.Ya. Shneiberg, Zeros of integrals along trajectories of ergodic systems, Funct. Anal. Appl., 19:2 (1985), 160-161;

И.Я. Шнейберг, Нули интегралов вдоль траекторий эргодических систем, Функц. анализ и его прил., 19:2 (1985), 92-93


\bibitem{De}
N. V. Denisova, Recurrence of integrals of conditionally periodic functions, Dokl. Math., 108:1 (2023), 316-319;

Н.В. Денисова, Возвращаемость интегралов условно периодических функций, Докл. РАН. Матем., информ., проц. упр., 512 (2023),  85-88

\bibitem{11} V.V. Ryzhikov, Recurrence of integral zeros on trajectories of ergodic flow. Preprint (2024)

\bibitem{Mo}Н.Г. Мощевитин, О возвращаемости интеграла гладкой трехчастотной условнопериодической функции, Матем. заметки, 58:5 (1995), 723-735; 

N.G. Moshchevitin, Recurrence of the integral of a smooth three-frequency conditionally periodic function, Math. Notes, 58:5 (1995), 1187-1196

\bibitem{Kon}С.В. Конягин, О возвращаемости интеграла нечетной условнопериодической функции, Матем. заметки, 61:4 (1997), 570-577; 

S.V. Konyagin, Recurrence of the integral of an odd conditionally periodic function, Math. Notes, 61:4 (1997), 473-479


\bibitem{Koz}	В.В. Козлов, Методы качественного анализа в динамике твердого тела, РХД, Ижевск, 2000 

V. V. Kozlov, Methods of qualitative analysis in the dynamics of a rigid body.
Scientific Publishing Center Regular and Chaotic Dynamics, Izhevsk, 2000. 248 pp. 

\bibitem{Ko}А.В. Кочергин, О росте сумм Биркгофа над поворотом окружности, Матем. заметки, 113:6 (2023), 836-848; 

A.V. Kochergin, On the Growth of Birkhoff Sums over a Rotation of the Circle, Math. Notes, 113:6 (2023), 784-793


\bibitem{97}В.\,В.~Рыжиков, Полиморфизмы, джойнинги и тензорная простота динамических систем, Функц. анализ и его прил., 31:2 (1997),  45-57

V.V. Ryzhikov, 
Polymorphisms, joinings, and the tensor simplicity of dynamical 
systems, Funct. Anal. Appl., 31:2 (1997), 109-118

\bibitem{23}
V.V. Ryzhikov, Self-joinings and generic extensions of ergodic systems, Funct. Anal. Appl., 57:3 (2023), 236-247 

В. В. Рыжиков, Самоприсоединения и типичные расширения эргодических систем, Функц. анализ и его прил., 57:3 (2023),  74-88 

\bibitem{We} B. Weiss, 
Single orbit dynamics. BMS Regional Conference Series in Mathematics 95. Providence, RI: AMS, 113 p. (2000).

\end{thebibliography}
\end{document}